\documentclass[journal]{IEEEtran}

\usepackage{cite}
\usepackage{amsmath,amssymb,amsfonts,bm}
\usepackage{bbm}
\usepackage{graphicx}
\usepackage{subfigure}
\usepackage{textcomp}
\usepackage{algorithm}
\usepackage{algpseudocode}

\makeatletter
\renewcommand*\env@matrix[1][*\c@MaxMatrixCols c]{%
	\hskip -\arraycolsep
	\let\@ifnextchar\new@ifnextchar
	\array{#1}}
\makeatother
\usepackage{mathrsfs}
\usepackage{xcolor}
\allowdisplaybreaks[4]

\newtheorem{theorem}{Theorem}
\newtheorem{lemma}{Lemma}
\newtheorem{assum}{Assumption}

\newenvironment{proof}{\begin{IEEEproof}}{\end{IEEEproof}}

\begin{document}
	\title{Policy Gradient Methods for the Cost-Constrained LQR: Strong Duality and Global Convergence
		\thanks{Research of F. Zhao and K. You was supported by National Key R\&D Program of China (2022ZD0116700) and National Natural Science Foundation of China (62033006, 62325305). (Corresponding author: Keyou You)}
		\thanks{F. Zhao and K. You are with the Department of Automation and Beijing National Research Center for Information Science and Technology, Tsinghua University, Beijing 100084, China. (e-mail: zhaofr18@mails.tsinghua.edu.cn, youky@tsinghua.edu.cn)}
		\author{Feiran Zhao, Keyou You}}
	\maketitle
\begin{abstract}
In safety-critical applications, reinforcement learning (RL) needs to consider safety constraints. However, theoretical understandings of constrained RL for continuous control are largely absent. As a case study, this paper presents a cost-constrained LQR formulation, where a number of LQR costs with user-defined penalty matrices are subject to constraints. To solve it, we propose a policy gradient primal-dual method to find an optimal state feedback gain. Despite the non-convexity of the cost-constrained LQR problem, we provide a constructive proof for strong duality and a geometric interpretation of an optimal multiplier set. By proving that the concave dual function is Lipschitz smooth, we further provide convergence guarantees for the PG primal-dual method. Finally, we perform simulations to validate our theoretical findings.
\end{abstract}

	
\section{Introduction}
Recent years have witnessed tremendous successes of reinforcement learning (RL) in continuous control and sequential decision-making tasks. In some safety-critical applications, e.g., autonomous driving \cite{fisac2018general}, robotics~\cite{ono2015chance}, and finance~\cite{krokhmal2002portfolio}, constraints must be taken into account as an additional learning objective. This leads to the conceptualization of constrained Markov decision processes (CMDPs)~\cite{altman1999constrained}, where a number of safety-related cost functions are subject to constraints.

Policy gradient (PG) method, as an essential approach of RL,  parameterizes the policy and directly updates it to minimize a cost function. Recently, PG methods have been used to solve CMDP problems~\cite{paternain2019constrained,chow2017risk, ding2020natural, pan2019risk,yu2019convergent}. In particular, the PG primal-dual approach has attracted increasing attentions due to its simplicity~\cite{paternain2019constrained,chow2017risk, ding2020natural}. It alternates between the primal iteration (minimizing the Lagrangian with PG descent) and the dual iteration (updating the multiplier with subgradient ascent) to find an optimal policy-multiplier pair. A key to the convergence of PG primal-dual methods is strong duality of CMDPs. However, strong duality has only been shown for the CMDPs with finite state-action space~\cite{altman1999constrained} or bounded cost functions~\cite{paternain2019constrained}, and it remains unclear for continuous control scenarios. This is because the duality analysis is notoriously challenging for non-convex continuous optimization problems. Consequently, the convergence guarantees of PG primal-dual methods for continuous control are largely absent.
% ¶ÔÅŒÐÔºÍÊÕÁ²ÐÔ»»Ò»ÏÂË³Ðò
%While the convergence of PG primal-dual methods has been shown for CMDPs with finite state-action space~\cite{altman1999constrained} or bounded cost functions~\cite{paternain2019constrained}, it remains unclear for the continuous control tasks. The bottleneck lies in that the duality analysis is notoriously challenging for non-convex optimization problems.

To improve theoretical understandings of PG methods for continuous control, there has been an increasing interest in studying their performance on classical control problems, e.g., the celebrated linear quadratic regulator (LQR) problem~\cite{fazel2018global,zhao2024data,zhang2021policy}. As a case study of continuous CMDPs, our previous work \cite{zhao2021global} considers the PG primal-dual method for the LQR with a single cost constraint. While strong duality is proved, the proof techniques in \cite{zhao2021global} cannot be applied to the multiple-constraint case. Moreover, the dual problem is a single variable optimization problem, the convergence analysis of which only provides limited insights for general CMDPs with multiple constraints.

In this paper, we propose a cost-constraint LQR formulation, where a number of LQR costs with user-defined penalty matrices are upper bounded. It can be viewed as a discrete-time counterpart of the LQR with integral quadratic constraints in \cite{lim1999stochastic}. Compared with \cite{zhao2021global}, this formulation allows multiple cost constraints and hence serves as an ideal benchmark for studying continuous CMDPs with unbounded cost. The cost-constraint LQR is also an instance of multi-objective control, which automatically balances different control objectives, e.g., mean performance and variance in risk-sensitive control~\cite{tsiamis2020risk,zhao2021infinitehorizon}. Applications include energy-constrained building control \cite{lee2018simulation} and the high-performance aircraft control \cite{kreisselmeier1983application}.

To solve the cost-constraint LQR problem, we propose a PG primal-dual method to find an optimal state feedback gain. The primal iteration uses PG methods to minimize the Lagrangian, which is a weighted LQR cost. While the convergence of the primal iteration is established in~\cite{fazel2018global}, there are two main challenges in analyzing the convergence of the dual iteration. First, we need to show strong duality for the non-convex cost-constraint LQR problem. Second, the analytic properties of the dual function (e.g., differentiability and smoothness) are largely unclear, limiting the attainable convergence rate of the PG primal-dual method. To this end, we first prove that the minimizer of the Lagrangian is unique and continuous in the multiplier, based on which we prove the strong duality. In particular, our proof is constructive, i.e., we construct a feasible optimization problem, whose solution is an optimal multiplier with zero duality gap. Then, we show that the dual function is differentiable, and the unique subgradient is actually the gradient. By using perturbation theory for algebraic Riccati equation, we further show that the dual function is Lipschitz smooth and provide convergence guarantees for the PG primal-dual method. As a comparison, the work \cite{zhao2021global} does not discover the smoothness of the dual function and hence only prove a slower convergence rate. We hope our work paves the way for rigorously understanding PG methods for general continuous CMDPs with unbounded cost functions.

The rest of this paper is organized as follows. Section II formulates the cost-constrained LQR problem. Section III proposes the PG primal-dual method and shows strong duality. Section IV shows its convergence. Section V performs simulations to validate the theoretical results. Conclusion is made in Section VI.

\textbf{Notations.} We use $I_n$ to denote the $n$-by-$n$ identity matrix. We use $\underline{\sigma}(\cdot)$ to denote the minimal singular value of a matrix. We use $\|\cdot\|$ to denote the $2$-norm of a vector or matrix. We use $\rho(\cdot)$ to denote the spectral radius of a square matrix. We use $[x]_+$ to denote the projection of $x\in\mathbb{R}^n$ onto $\mathbb{R}^n_+$.

\section{Problem formulation}

Consider the following discrete-time linear system
\begin{equation}\label{equ:sys}
x_{t+1} = Ax_t + Bu_t, ~~x_0 \sim \mathcal{D},
\end{equation}
where $x_t \in \mathbb{R}^n$ is the state, $u_t\in \mathbb{R}^m$ is the control input. The pair $(A,B)$ is controllable. The initial state $x_0$ is sampled from a distribution $\mathcal{D}$ with zero mean and $\mathbb{E}[x_0x_0^{\top}]=1$.

We consider the cost-constrained LQR as a case study of continuous constrained Markov decision processes (CMDP) with unbounded cost functions~\cite{altman1999constrained}. Specifically, we aim to find an optimal policy sequence $\{\pi_t\}$ that solves
\begin{equation}\label{prob:dis_clqr}
\begin{aligned}
&\mathop{\text{min}}\limits_{\{\pi_t\}} ~J_0 := \mathbb{E}_{x_0} \left[ \sum_{t=0}^{\infty}(x_{t}^{\top} Q_0 x_{t}+u_{t}^{\top} R_0 u_{t}) \right]\\
&~\text {s.t.} ~ J_i:= \mathbb{E}_{x_0} \left[ \sum_{t=0}^{\infty}(x_{t}^{\top} Q_i x_{t}+u_{t}^{\top} R_i u_{t})\right] \leq c_i,  \\
&~~~u_t = \pi_t(x_0, u_0,\cdots, x_t),~\text{and}~(\ref{equ:sys}), ~i \in \{1,2,\dots,N\},
\end{aligned}
\end{equation}
where  $Q_i, R_i, i\in  \{0, 1,\dots,N\}$ are user-defined penalty matrices satisfying the following assumption.
\begin{assum}
	\label{assumption}
	The penalty matrices satisfy $Q_0, R_0 > 0$ and $Q_i, R_i \geq 0, \forall i \in \{1,2,\dots,N\}$.
\end{assum}
 
This formulation can also be viewed as an instance of multi-objective  control~\cite{chen2021multi,lim1999stochastic}, which automatically balances different control objectives without manual weight-tuning. It can be shown by using the techniques in \cite{zhao2021infinitehorizon} that the optimal policy of \eqref{prob:dis_clqr} is linear state feedback, i.e., $u_t = K^* x_t$. Thus, we use $K\in \mathbb{R}^{m\times n}$ to parameterize the policy and focus on the following optimization problem to find $K^*$:
\begin{equation}\label{prob:clqr}
\mathop{\text{min}}\limits_{K} ~J_0(K), 
~\text{s.t.} ~ J_i(K) \leq c_i,  ~i \in \{1,2,\dots,N\}.
\end{equation}

%In particular, it covers the risk-constrained LQR \cite{tsiamis2020risk,zhao2021infinitehorizon} balancing mean performance and variance, the energy-constrained LQR in building control~\cite{lee2018simulation}, and the high-performance aircraft control in different flying conditions \cite{kreisselmeier1983application} .

%It automatically balances conflicting control objectives~\cite{chen2021multi}, such as mean performance and variance in risk-sensitive control~\cite{tsiamis2020risk}, regulation and consensus in distributed control~\cite{li2019distributed}. 

%Moreover, the cost-constrained LQR can be viewed as a multi-objective $\mathbb{H}_2$ control problem~\cite{chen2021multi}. 

Since $J_i(K)$ is non-convex~\cite{fazel2018global}, the problem \eqref{prob:clqr} is a constrained non-convex problem. While the convergence of policy gradient (PG) methods has been shown for the LQR~\cite{fazel2018global}, the non-convex constraints render the analysis of \eqref{prob:clqr} more involved. In this paper, we leverage primal-dual optimization theory to show the convergence of PG methods for \eqref{prob:clqr}.

\section{Policy gradient primal-dual methods for the cost-constrained LQR}
In this section, we first propose the policy gradient primal-dual method to solve the cost-constrained LQR problem.  Then, we show strong duality by analyzing properties of the Lagrangian.
\subsection{The policy gradient primal-dual method to solve \eqref{prob:clqr}}
Define the stabilizing set 
$
\mathcal{S} = \{K \in \mathbb{R}^{m\times n}|\rho(A-BK)<1\}.
$
For $K\in \mathcal{S}$, the LQR costs are finite, i.e., $J_i(K)<\infty, \forall i \in \{0,1,\dots,N\}$.

Let $\lambda = [\lambda_1,\lambda_2,\dots,\lambda_N]^{\top} \geq 0$ be the multiplier of the cost-constrained LQR (\ref{prob:clqr}). Then, the Lagrangian is given by
\begin{equation}\label{def:L}
\begin{aligned}
\mathcal{L}(K,\lambda) &:= J_0(K)+ \sum_{i=1}^{N}\lambda_i(J_i(K) -c_i )\\
&=\mathbb{E} \left[ \sum_{t=0}^{\infty}(x_{t}^{\top} Q_\lambda x_{t}+u_{t}^{\top} R_\lambda u_{t}) \right] - \sum_{i=1}^{N}\lambda_ic_i,
\end{aligned}
\end{equation}
where  $Q_{\lambda} := Q_0+ \sum_{i=1}^{N} \lambda_iQ_i $ and $R_{\lambda} := R_0+ \sum_{i=1}^{N} \lambda_iR_i $ are weighted penalty matrices. The dual function is given by
\begin{equation}\label{def:D}
D(\lambda) =  \min \limits_{K \in \mathcal{S}} \mathcal{L}(K,\lambda).
\end{equation}

Our PG primal-dual method alternates between
\begin{subequations}\label{equ:primal-dual}
	\begin{align}
	&K^{k}\in \mathop{\text{argmin}} \limits_{K \in \mathcal{S}}~ \mathcal{L}(K,\lambda^k), \label{prim_iterate} \\
	&\lambda^{k+1}= [\lambda^k+ \eta^k \cdot d^k]_+. \label{dual_iterate}
	\end{align}
\end{subequations}
The primal iteration \eqref{prim_iterate} is solved by PG methods, i.e., for a fixed multiplier $\lambda^k$ we iterate
\begin{equation}\label{equ:PG}
K^{+} = K - \zeta \nabla_K \mathcal{L}(K,\lambda^k),
\end{equation}
where $\zeta>0$ is a stepsize. The dual iteration in \eqref{dual_iterate} uses projected subgradient ascent to update the multiplier, where $\eta^k>0$ is a stepsize, and $d^k$ is a subgradient of $D(\lambda)$ at $\lambda^k$ computed using $K^k$.  

While the convergence of PG methods \eqref{equ:PG} for the primal iteration has been shown in \cite{fazel2018global}, the optimization landscape of the dual iteration \eqref{dual_iterate} remains unclear. A key to establishing the convergence of \eqref{dual_iterate} is the strong duality between the primal problem (\ref{prob:clqr}) and the following dual problem 
\begin{equation}\label{prob:dual}
\mathop{\text{maximize}}\limits_{\lambda \geq 0}~ D(\lambda).
\end{equation}
However, since the problem \eqref{prob:clqr} is non-convex, strong duality does not trivially hold. In the sequel, we show the strong duality by proving analytic properties of the Lagrangian.

%Since (\ref{prob:new_rclqr}) is non-convex, it does not trivially hold. {Even though the Lagrangian in \eqref{def:L} is the LQR cost with a linear term,  its non-convex optimization landscape is yet unclear.} Thus, computing the primal update in \eqref{prim_iterate} is itself challenging.  In the rest of this section, we show that (a) $\mathcal{L}(X,\mu)$ is coercive over $\mathcal{S}$ and locally gradient dominated in Section \ref{subsec:property}, which is key to establish that a critical point of \eqref{prim_iterate} is globally optimal; (b) $\mathcal{L}(X,\mu)$ and its gradient are locally Lipschitz in Section \ref{subsec_local}, which implies a linear convergence rate of gradient methods for solving \eqref{prim_iterate}; (c) The strong duality property indeed holds in Section \ref{subsec_strong}. Combining these results proves the global convergence of \eqref{primal-dual}. Note that all the proofs on the properties of the Lagrangian are provided in Appendix \ref{Appendix:A} and \ref{apx:lip}. 

\subsection{Strong duality between the primal problem (\ref{prob:clqr}) and the dual problem \eqref{prob:dual}}	
	
The Lagrangian $\mathcal{L}(K,\lambda)$ in \eqref{def:L} is a standard LQR cost with a constant bias. Thus, for a fixed multiplier $\lambda$ the minimizer of $\mathcal{L}(K,\lambda)$ is uniquely given by
\begin{equation}\label{def:K}
K_{\lambda}^*  = (R_{\lambda} + B^{\top}P_{\lambda}^*B )^{-1}BP_{\lambda}^*A,
\end{equation}
where $P_{\lambda}^*$ is the positive definite solution to the algebraic Riccati equation (ARE)
\begin{equation}\label{equ:dare}
P_{\lambda}^*=A^{\top} P_{\lambda}^* A+Q_\lambda-A^{\top} P_{\lambda}^* B(R_{\lambda}+B^{\top} P_{\lambda}^* B)^{-1} B^{\top} P_{\lambda}^* A.
\end{equation}

We first show the continuity of both $K_{\lambda}^*$ and  $J_i(K_{\lambda}^*)$ in $\lambda$, the proof of which is provided in Appendix \ref{app:1}.
\begin{lemma}\label{lem:continuity}
	The unique minimizer of the Lagrangian $K_{\lambda}^*$ and the constrained costs $J_i(K_{\lambda}^*), \forall i \in \{1,2,\cdots, N\}$ are continuous in $\lambda$ over $\mathbb{R}^N_+$.
\end{lemma}

The continuity in Lemma \ref{lem:continuity} is a strong result and usually does not hold for general constrained Markov decision processes (CMDPs). In particular, it requires the uniqueness of the minimizer of the Lagrangian as a necessary condition. 

Our proof of strong duality relies on the construction of the following optimization problem
\begin{equation}\label{equ:pareto}
\begin{aligned}
&\mathop{\text{minimize}}\limits_{\lambda}~~ z^{\top}\lambda, \\
&\text{subject to}~~ J_i(K_{\lambda}^*) \leq c_i, ~i\in \{1,\dots, N\},
\end{aligned}
\end{equation}
where $z>0$ is a constant vector. We have the following result regarding \eqref{equ:pareto}.
\begin{lemma}\label{lem:pareto}
	Suppose that Slater's condition holds for the cost-constrained LQR problem \eqref{prob:clqr}, i.e., there exists a policy $\widetilde{K} \in \mathcal{S}$ such that $J_i(\widetilde{K})<c_i, \forall i \in \{1,2,\cdots, N\}$. Then, the problem \eqref{equ:pareto} is feasible and an optimal solution exists.
\end{lemma}
\begin{proof}
First, we use contradiction to prove the feasibility of \eqref{equ:pareto}, i.e., there exists a multiplier $\lambda$ such that $J_i(K_{\lambda}^*) \leq c_i, ~\forall i\in \{1,\dots, N\}$. Without loss of generality, suppose that  for any $\lambda \geq 0$ we have $J_i(K_{\lambda}^*) > c_i, \forall i \in \{1,\dots,j\}$. Then, we let $\lambda_i = 0, \forall i \in \{j+1, \dots, N \}$. By Slater's condition, there exists some constants $a_i>0$ such that $J_i(\widetilde{K})+a_i \leq c_i, \forall i \in \{1,2,\dots, N\}$. Then, the following relations hold
\begin{align*}
J(\widetilde{K}) &\geq D(\lambda) - \sum_{i=1}^{N} \lambda_i (J_i(\widetilde{K}) - c_i)\\
&\geq J(K_{\lambda}^*) + \sum_{i=1}^{N} \lambda_i (J_i(K_{\lambda}^*) - J_i(\widetilde{K})) \\
&\geq J(K_{\lambda}^*) + \sum_{i=1}^{N} \lambda_i (J_i(K_{\lambda}^*) - c_i + a_i) \\
&> J(K_{\lambda}^*) + \sum_{i=1}^{k}\lambda_ia_i,
\end{align*}
where the first two inequalities hold by the definition of $D(\lambda)$ in (\ref{def:D}), the third inequality follows from Slater's condition and the fourth one from the hypothesis and the fact $\lambda_i = 0, \forall i \in \{j+1, \dots, N \}$. Letting $\lambda_i \rightarrow \infty, \forall i \in \{1,2,\dots,k \}$ yields $J(\widetilde{K}) > \infty$, which contradicts the Slater's condition. Thus, there must exist a multiplier $\widetilde{\lambda}$ such that $J_i(K(\widetilde{\lambda})) \leq c_i, ~\forall i\in \{1,\dots, N\}$. 

Next, we show that an optimal solution to \eqref{equ:pareto}  exists. Consider the compact set $\mathcal{C} = \{\lambda \geq 0| \lambda \leq \widetilde{\lambda} \}$. Clearly, the solution to the problem (\ref{equ:pareto}) must lie in $\mathcal{C}$. Thus, we only need to focus on the compact set $\mathcal{C}$. By Lemma \ref{lem:continuity}, both the objective and constraints in problem (\ref{equ:pareto}) are continuous over $\mathcal{C}$. Then, by Weierstrass' Theorem \cite[Proposition A.8]{bertsekas1997nonlinear}, the set of minima of the problem (\ref{equ:pareto}) over $\mathcal{C}$ is nonempty, which implies the existence of an optimal solution.
\end{proof}

Based on Lemmas \ref{lem:continuity} and \ref{lem:pareto}, we show that strong duality holds, and an optimal multiplier to the dual problem is given exactly by an solution to \eqref{equ:pareto}.

\begin{theorem}[\textbf{Strong duality}]\label{theorem:duality}
	Under Slater's condition, there is no duality gap between the primal problem (\ref{prob:clqr}) and the dual problem \eqref{prob:dual}. Moreover, an optimal-multiplier pair is given by $(K_{\lambda^*}^*, \lambda^*)$, where $\lambda^*$ is an optimal solution to \eqref{equ:pareto}.
\end{theorem}
\begin{proof}
	Let $
	\lambda^* = [\lambda_1^*, \dots,\lambda_N^*]^{\top}
	$
 	be an optimal solution to \eqref{equ:pareto}, the existence of which is guaranteed by Lemma \ref{lem:pareto}. We show that the policy-multiplier pair $(K^*, \lambda^*)$ with $K^*:=K_{\lambda^*}^*$ satisfies the optimality conditions~\cite[Proposition 6.1.5]{bertsekas1997nonlinear}
	\begin{subequations}\label{equ:saddle}
	\begin{align}
	&\mathcal{L}(K^*,\lambda^*) = \min_{K \in \mathcal{S}}\mathcal{L}(K,\lambda^*), \label{equ:first} \\
	&J_i(K^*) \leq c_i, \label{equ:second}  \\
	&\lambda_i^*(J_i(K^*)- c_i)=0, ~\forall i \in \{1,2,\cdots, N\}, \label{equ:third} 
	\end{align}
	\end{subequations}
	which directly implies the strong duality~\cite[Chapter 6]{bertsekas1997nonlinear}.
	
	By the definition of $K^*$ and $\lambda^*$, the first two optimality conditions (\ref{equ:first})-(\ref{equ:second}) hold, and it remains to show \eqref{equ:third}. We first prove that the constrained cost function $J_i(K_{\lambda}^*)$ is monotone non-increasing in $\lambda_i$. Define 
	$\lambda = [\lambda_1, \dots, \lambda_j,\dots,\lambda_N]^{\top}$ and $\lambda' = [\lambda_1, \dots, \lambda_j',\dots,\lambda_N]^{\top}$. Then, it follows from the definition of the dual function in \eqref{def:D} that
	\begin{align*}
	&J(K_{\lambda}^*) + \sum_{i=1}^{N}\lambda_{i}J_i(K_{\lambda}^*) \leq 
	J(K_{\lambda'}^*) + \sum_{i=1}^{N}\lambda_{i}J_i(K_{\lambda'}^*), \\
	&J(K_{\lambda}^*) + \sum_{i=1}^{N}\lambda_{i}J_i(K_{\lambda}^*) -\lambda_jJ_j(K_{\lambda}^*) + \lambda_j'J_j(K_{\lambda}^*)  \\
	&\geq
	J(K_{\lambda'}^*) + \sum_{i=1}^{N}\lambda_{i}J_i(K_{\lambda'}^*) -\lambda_jJ_j(K_{\lambda'}^*) + \lambda_j'J_j(K_{\lambda'}^*).
	\end{align*}
	Then, subtracting in both sides yields that
	$$
	(\lambda_j' - \lambda_j)(J_j(K_{\lambda}^*) - J_j(K_{\lambda'}^*) ) \geq 0.
	$$
	
	Now, we show that \eqref{equ:third} holds. For $i$-th constraint, consider two cases. If $\lambda_i^* = 0$, then it trivially holds.  If $\lambda_i^* >0$, then it follows from the definition of $\lambda^*$ that $$J_i\left(K\left([\lambda_1^*,\cdots,\lambda_{i-1}^*, 0,\lambda_{i+1}^*,\cdots,\lambda_N^*]^{\top}\right)\right) > c_i.$$ 
	Since $J_i(K_{\lambda}^*)$ is continuous and monotone non-increasing in $\lambda_i$, we have $J_i(K_{\lambda^*}^*)=c_i$ and the $i$-th equality in \eqref{equ:third} holds. Since the above derivation holds for $i \in \{1,2,\cdots, N\}$, the third optimality condition \eqref{equ:third} holds. By \cite[Proposition 6.1.5]{bertsekas1997nonlinear}, strong duality holds, and $(K_{\lambda^*}^*, \lambda^*)$ is an optimal policy-multiplier pair.
\end{proof}

\begin{figure}[t]
	\centerline{\includegraphics[width=50mm]{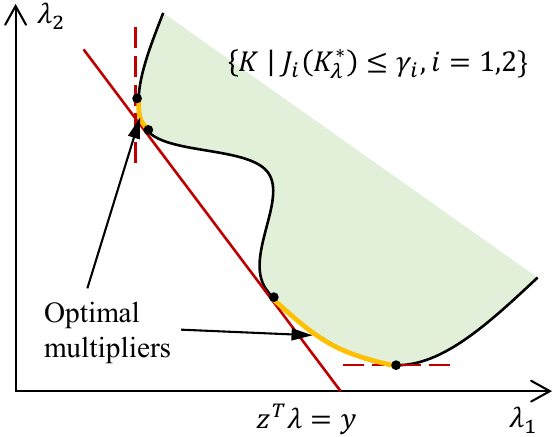}}
	\caption{A geometry interpretation for optimal multipliers in a two-constraint example, where $y$ is some positive constant.}
	\label{pic:pareto}
\end{figure}
We make some remarks on Theorem \ref{theorem:duality}. First, compared with the strong duality result for the LQR with a single cost constraint in \cite[Theorem 2]{zhao2021global}, Theorem \ref{theorem:duality} is more general as it considers multiple cost constraints. Second, the proof of Theorem \ref{theorem:duality} is constructive, i.e., we show that an optimal multiplier is given by the solution to the optimization problem \eqref{equ:pareto}. Notice that since the parameter $z >0$ of \eqref{equ:pareto} is arbitrarily given, an optimal multiplier may not be unique. In fact, the solutions to \eqref{equ:pareto} for all $z >0$ belong to the optimal multiplier set, and a geometry interpretation is provided in Fig. \ref{pic:pareto}. Third, Theorem \ref{theorem:duality} has independent interests of its own. In particular, the proof techniques can be applied to show strong duality of more general constrained optimization problems, as long as the continuity properties (c.f. Lemma \ref{lem:continuity}) hold.

However, the non-convex problem \eqref{equ:pareto} is intractable, and hence cannot be used to obtain an optimal multiplier. In the sequel, we use the PG primal-dual method in \eqref{equ:primal-dual} to find an optimal policy-multiplier pair.

\section{Convergence of the PG primal-dual method for the cost-constrained LQR}
In this section, we first provide analytic properties of the dual function, i.e., differentiability and Lipschitz smoothness. Then, we show the convergence of the PG primal-dual method in \eqref{equ:primal-dual} for the cost-constrained LQR.
\subsection{Properties of the dual function}
First, by convex optimization theory~\cite[Chapter 5]{bertsekas1997nonlinear}, the dual function $D(\lambda)$ is concave  over $\lambda \in \mathbb{R}_+^N$. Second, since the minimizer of the Lagrangian function is unique, the dual function is  differentiable.
\begin{lemma}[\textbf{Differentiability of the dual function}]\label{lem:gradient}
	The dual function $D(\lambda)$ is differentiable over $\mathbb{R}^N_+$, and its gradient is 
	\begin{equation}\label{equ:subgrad}
	\nabla D(\lambda) =  \begin{bmatrix}
	J_1(K_{\lambda}^*)-c_1 &  \cdots & J_N(K_{\lambda}^*)-c_N
	\end{bmatrix}^{\top}.
	\end{equation}
\end{lemma}
\begin{proof}
Since the Lagrangian function has a unique minimizer \eqref{def:K}, the subdifferential of the concave dual function $D(\lambda)$ is a singleton, and $D(\lambda)$ is differential over $\lambda \in \mathbb{R}_+^N$. Moreover, the gradient (actually the unique subgradient) can be computed by Danskin theorem \cite[Chapter 5]{bertsekas1997nonlinear}.
\end{proof}
  
%The differentiability of the dual function in Lemma \ref{lem:gradient} is a strong property that usually lacks in the CMDP problems, where only a subgradient is accessible. 

By Lemma \ref{lem:gradient}, the dual iteration \eqref{dual_iterate} is projected gradient ascent. Next, we further show the Lipschitz continuous of the gradient $\nabla D(\lambda)$ describing the smoothness of $D(\lambda)$. By Theorem \ref{theorem:duality}, an optimal multiplier $\lambda^*$ is finite, i.e., there exists a compact set $\Omega$ such that $\lambda^* \in \Omega$. Thus, it suffices to focus on $\Omega$ for the smoothness analysis.

\begin{lemma}[\textbf{Local Lipschitz smoothness}]\label{lem:sublipschitz}
	For any $\lambda, \lambda'\in \Omega$, there exist constants $\gamma_1>0$ and $\mu >0$ such that if $\|\lambda'-\lambda\| \leq \gamma_1$, 
	then it holds that 
	$
	\|\nabla D(\lambda') - \nabla D(\lambda)\| \leq \mu \|\lambda' -{\lambda}\|.
	$
\end{lemma}

The proof of Lemma \ref{lem:sublipschitz} is provided in Appendix \ref{app:2}, which leverages the perturbation theory for Lyapunov equations~\cite{fazel2018global} and Riccati equations~\cite{sun1998perturbation}. As a key analytic property, the smoothness enables us to improve the convergence rate of the primal-dual method over that of \cite{zhao2021global}, where such a property is not proved. 

Based on Lemmas \ref{lem:gradient} and \ref{lem:sublipschitz}, we next show the convergence of the PG primal-dual method in \eqref{equ:primal-dual}.

\subsection{Convergence of the PG primal-dual method}\label{subsec:conver}
Since the Lagrangian is a weighted LQR cost, for a fixed multiplier the PG method \eqref{equ:PG} meets global convergence~\cite{fazel2018global}, i.e.,  the primal iteration \eqref{prim_iterate}  returns a solution lying in a compact set \cite{fazel2018global}
\begin{equation}\label{def:epsilon}
	K^k \in \{K| \|K-K^*_{\lambda^k} \|\leq  \epsilon\}
\end{equation}
for some constant $\epsilon>0$. Consequently, we can use only approximated gradient of the dual function in the dual iteration \eqref{dual_iterate}, i.e., 
$
d^k = \begin{bmatrix}
J_1(K^k)-c_1 &  \cdots & J_N(K^k)-c_N
\end{bmatrix}^{\top}.
$
To ensure the boundedness of the multiplier, the dual iteration \eqref{dual_iterate} projects on the compact set $\Omega$\footnote{There are manifold approaches to finding $\Omega$; see \cite{nedic2009subgradient,ding2020natural}. We omit it here and focus on the convergence analysis due to space limitation.} containing $\lambda^*$, i.e., 
\begin{equation}\label{equ:appro_dual}
\lambda^{k+1}= \Pi_{\Omega} (\lambda^k+ \eta^k   d^k).
\end{equation} 

We first provide bounds for the norm of $d^k$ and the distance between $d^k$ and $\nabla D(\lambda^k)$, the proof of which is provided in Appendix \ref{app:3}.
\begin{lemma}\label{lem:sub_error}
For all $k \in \mathbb{N}$, there exist uniform constants $\gamma_2>0$, $\bar{d}>0$, and $c>0$ such that if $\epsilon \leq  \gamma_2$, then $\|d^k - \nabla D(\lambda^k)\| \leq c\epsilon$ and $\|d^k\|\leq \bar{d}$. 
\end{lemma}

Since $\Omega$ is a compact set, we let $\omega = \max_{\lambda \in \Omega} \{\|\lambda\|\}$. Define the regret of the dual function for $k \geq 1$ as
$$
\text{Regret}_k :=\frac{1}{k}\sum_{i=0}^{k-1}(D^*-D(\lambda^{i+1})).
$$
Then,  we show the convergence of the primal-dual method.

\begin{theorem}[\textbf{Global convergence}]\label{thm:finald}
If $\epsilon \leq  \gamma_2$ and $0<\eta \leq \min \{\gamma_1/\bar{d}, 2/\mu \}$, then for $k \in \mathbb{N}_+$ it holds that
$$
\text{Regret}_k \leq \frac{\sqrt{p_2 D^*}}{\sqrt{k}} + \sqrt{p_1p_2}.
$$
where  $p_1 =  \eta c \epsilon(1-{\mu \eta}/2)(\eta c \epsilon + 2\bar{d} ) + c\epsilon \omega$ and $p_2 = (1-{\mu \eta}/{2}){\omega^2}/{\eta^2}$.
\end{theorem}

The proof is provided in Appendix \ref{app:3}.
Theorem \ref{thm:finald} shows that when the solution error $\epsilon$ of the primal iteration is sufficiently small, the regret is upper bounded by two terms signifying a sublinear decrease and a bias  polynomial in $\epsilon$. Moreover, the sublinear decrease matches the attainable rate $\mathcal{O}(1/\sqrt
{k})$ of first-order methods in online convex optimization of smooth functions~\cite[Chapter 3]{hazan2016introduction}. Compared with \cite[Theorem 6]{zhao2021global}, Theorem \ref{thm:finald} explicitly characterizes the effects of $\epsilon$ and shows a faster convergence rate thanks to the smoothness property in Lemma \ref{lem:sublipschitz}.

\section{Simulations}
This section validates the convergence of the PG primal-dual method for the cost-constrained LQR \eqref{prob:clqr}.
\subsection{Simulation model}
Consider an unmanned aerial vehicle operating on a 2-D plane. The dynamical model is given by a double integrator  
\begin{equation}\label{def:model}
A=\begin{bmatrix}
1 & 0.5 & 0 & 0 \\
0 & 1 & 0 & 0 \\
0 & 0 & 1 & 0.5 \\
0 & 0 & 0 & 1
\end{bmatrix}, B = \begin{bmatrix}
0.125 & 0 \\
0.5 & 0 \\
0 & 0.125 \\
0 & 0.5
\end{bmatrix},
\end{equation}
where $(x_{t,1}, x_{t,3})$ denote the position, $(x_{t,2}, x_{t,4})$ denote the velocity, and $u_k$ represents the acceleration. 

The task is to regulate the position to zero, while keeping the positions in two directions close at a low energy cost. To this end, for the cost-constrained LQR problem \eqref{prob:clqr} we select $c_1 = 10$, $c_2 =6$, $Q_0 = I_4, R_0 = I_2, Q_1 = 0, R_1 = 2\times I_2, R_2 = 0$, and
$$
Q_2 = \begin{bmatrix}
1 & 0 & -1 & 0\\
0 & 0 & 0 & 0 \\
-1 & 0 & 1 & 0\\
0 & 0 & 0 & 0
\end{bmatrix},
$$
where the first constraint is to limit the accumulated energy consumption, and the second is to make $x_{t,1}$ and $x_{t,3}$ close. 

In the sequel, we use the PG primal-dual method in \eqref{equ:primal-dual} to solve the cost-constrained LQR problem and validate its convergence.
\subsection{Convergence of the PG primal-dual method}
For the primal iteration \eqref{prim_iterate}, we select
$$K^0 = \begin{bmatrix}
0.5 & 0.5 & 0 & 0 \\
0 & 0 & 0.5 & 0.5
\end{bmatrix}
$$
as the initial stabilizing policy. We fix the constant stepsize to $\zeta = 10^{-3}$ in the PG update \eqref{equ:PG} and the number of PG steps per multiplier to $50,100,1000$, which approximately corresponds to the solution error $\epsilon \in \{10^{-1}, 10^{-2}, 10^{-8}\} $ in \eqref{def:epsilon}. For the dual iteration \eqref{equ:appro_dual}, we fix the constant stepsize to $\eta = 0.5$. 

\begin{figure}[t]
	\centerline{\includegraphics[width=55mm]{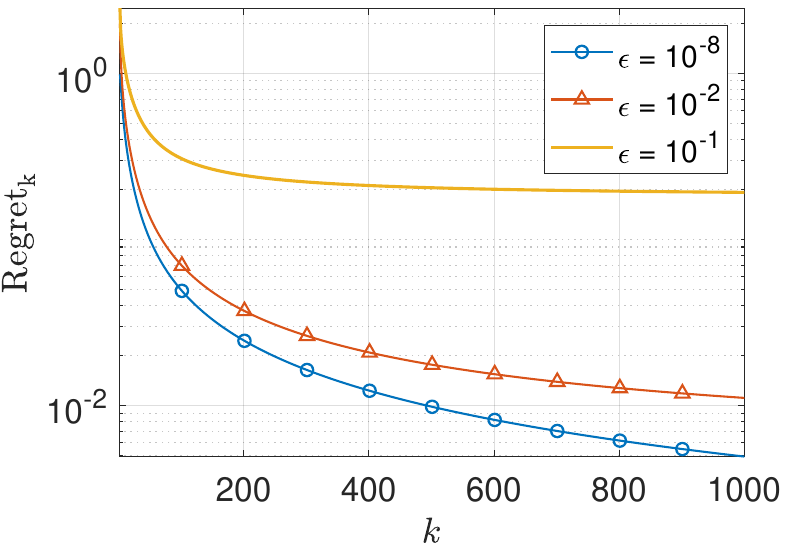}}
	\caption{Convergence of the dual iteration in the PG primal-dual method.}
	\label{pic:regret}
\end{figure}
Fig. \ref{pic:regret} illustrates the convergence of the dual iteration \eqref{equ:appro_dual}. As expected by Theorem \ref{thm:finald}, the regret converges at a sublinear rate, and the bias grows proportionally to the solution error $\epsilon$ in the primal iteration.

\begin{figure}[t]
	\centering
	\subfigure[Optimality gap.]{
		\includegraphics[width=50mm]{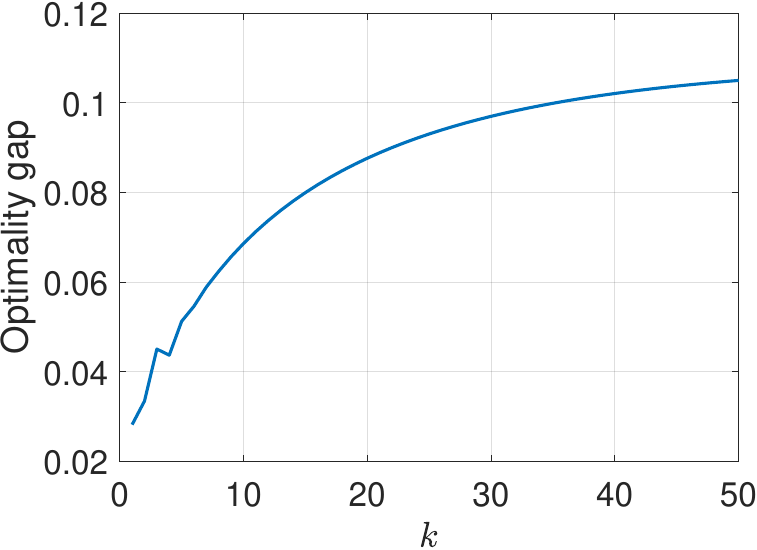}}
	
	\subfigure[Constraint violation.]{
		\includegraphics[width=50mm]{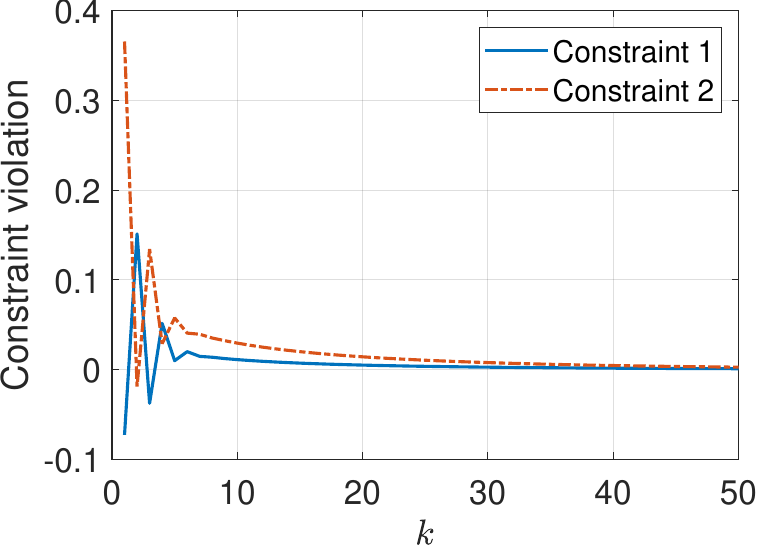}}
	\caption{Optimality gap and constraint violation of the PG primal-dual method.}
	\label{pic:pd}
\end{figure}

Next, we fix the number of PG steps per multiplier to $100$ and validate
the convergence of the optimality gap $(J_0(K^k)-J_0^*)/J_0^*$ and the constraint violation $(J_i(K^k)-c_i)/c_i, i\in\{1,2\}$, where $J_0^* = \min_{K \in \mathcal{S}} J_0(K)$. Fig. \ref{pic:pd} illustrates that in $50$ iterations, both two cost constraints are approximately satisfied, while the optimality gap slightly increases.

\section{Conclusion}
In this paper, we have proposed the cost-constrained LQR formulation, where the constraints are the LQR costs with user-defined penalty matrices. Then, we have proposed a policy gradient primal-dual method to solve the cost-constrained LQR problem. Despite the non-convexity of the primal-dual problem, we have shown strong duality and the Lipschitz smoothness of the dual function, based on which convergence guarantees have been provided. Finally, simulations have been performed to validate the theoretical results.

While the presented primal-dual method is model-based, it is straightforward to extend to the sample-based setting, where the gradient of the dual function is approximated using observations of the cost functions. We leave this interesting extension as an important future work.

\appendices
\section{Proof of Lemma \ref{lem:continuity}}\label{app:1}
	By the matrix inverse lemma, the ARE in (\ref{equ:dare}) is written as
\begin{equation}\label{equ:new_are}
P_{\lambda}^* = Q_{\lambda} + A^{\top} \left(I + P_{\lambda}^*G_{\lambda} \right)^{-1}P_{\lambda}^*A.
\end{equation}
We first show that $P_{\lambda}^*$ is a continuous function in $\lambda$ for $\lambda \geq 0$. To this end, we use implicit function theorem~\cite[Proposition A.25]{bertsekas1997nonlinear} to show that $\text{vec}(P_{\lambda}^*)$ is continuous in $\lambda$.

Vectorizing both sides of (\ref{equ:new_are}) yields that
\begin{align*}
&\text{vec}(P_{\lambda}^*)=\text{vec}(Q_{\lambda}) + \text{vec}(A^{\top}(I+ P_{\lambda}^*G_{\lambda} )^{-1}P_{\lambda}^*A )\\
&=\text{vec}(Q_{\lambda}) + (A^{\top} \otimes A) \text{vec}((I+ P_{\lambda}^*G_{\lambda} )^{-1}P_{\lambda}^*) \\
&=: \Phi_1(\text{vec}(P_{\lambda}^*), \lambda),
\end{align*}
where we defined the function $\Phi_1: \mathbb{R}^{n^2} \times \mathbb{R}^n \rightarrow \mathbb{R}^{n^2}$ and used the following relation
$
\text{vec}(ABC) = (C^{\top}\otimes A)\text{vec}(B)
$
for some matrices $A,B,C$ with proper dimensions. To apply the implicit function theorem, we first show that $\Phi_1(\text{vec}(P_{\lambda}^*), \lambda) - \text{vec}(P_{\lambda}^*)$ is continuous in both $P_{\lambda}^*$ and $\lambda$. This is clear since $Q_{\lambda}$ and $R_{\lambda}$ are linear functions of $\lambda$. Then, it remains to show 
\begin{equation}\label{equ:Phi}
\frac{\partial \Phi_1(\text{vec}(P_{\lambda}^*), \lambda) - \text{vec}(P_{\lambda}^*)}{\partial \text{vec}(P_{\lambda}^*)}
\end{equation}
is invertible. By \cite[(B.17)]{zhang2019policy}, it can be written as 
$$
\frac{\partial \Phi_1(\text{vec}(P_{\lambda}^*), \lambda) - \text{vec}(P_{\lambda}^*)}{\partial \text{vec}(P_{\lambda}^*)}=(A-BK_{\lambda}^*)\otimes(A-BK_{\lambda}^*)-I.
$$
Since $A-BK_{\lambda}^*$ is stable, the eigenvalues of $(A-BK_{\lambda}^*)\otimes(A-BK_{\lambda}^*)$ have absolute values smaller than one. Hence, (\ref{equ:Phi}) is invertible. Then, by the implicit function theorem, $P_{\lambda}^*$ is a continuous function of $\lambda$ over $\lambda \geq 0$. Thus, by the definition \eqref{def:K}, $K_{\lambda}^*$ is continuous in $\lambda$.

Next, we show that $J_i(K_{\lambda}^*), \forall i \in \{1,\dots,N\}$ are continuous in $\lambda$. By \cite{fazel2018global}, the LQR cost $J_i(K_{\lambda}^*)$ can be written as
\begin{equation}
J_i(K_{\lambda}^*) = \text{Tr}(P_i(K_{\lambda}^*))
\end{equation}
with $P_i(K_{\lambda}^*)$ the positive semi-definite solution to the Lyapunov equation
\begin{equation}\label{equ:P_i}
P_i(K_{\lambda}^*) = Q_i + K_{\lambda}^{*\top}R_i K_{\lambda}^* + (A-BK_{\lambda}^*)^{\top}P_i(K_{\lambda}^*)(A-BK_{\lambda}^*).
\end{equation}
Thus, it suffices to show that $P_i(K_{\lambda}^*)$ is continuous in $K_{\lambda}^*$. Vectorizing both sides of (\ref{equ:P_i}) yields that
\begin{align*}
&\text{vec}(P_i(K_{\lambda}^*)) = \text{vec}(Q_i) + \left(K_{\lambda}^{*\top} \otimes K_{\lambda}^*\right)\text{vec}(P_i(K_{\lambda}^*))\\
& + \left((A-BK_{\lambda}^*)^{\top} \otimes (A-BK_{\lambda}^*)\right)\text{vec}(P_i(K_{\lambda}^*))\\
&:=\Phi_2(\text{vec}(P_i(K_{\lambda}^*)), K_{\lambda}^*),
\end{align*}
where we defined the function $\Phi_2: \mathbb{R}^{n^2} \times \mathbb{R}^{m\times n} \rightarrow \mathbb{R}^{n^2}$. Clearly, $\Phi_2(\text{vec}(P_i(K_{\lambda}^*)), K_{\lambda}^*)-\text{vec}(P_i(K_{\lambda}^*)) $ is continuous in both $\text{vec}(P_i(K_{\lambda}^*))$ and $K_{\lambda}^*$. Analogously to \eqref{equ:Phi}, we have $
({\partial \Phi_2(\text{vec}(P_i(K_{\lambda}^*)) - \text{vec}(P_i(K_{\lambda}^*)) })/{\partial \text{vec}(P_i(K_{\lambda}^*))} 
$
is invertible. Thus, by the implicit function theorem, $\text{vec}(P_i(K_{\lambda}^*))$ is continuous in $K_{\lambda}^*$. Furthermore, $J_i(K_{\lambda}^*), \forall i \in \{1,\dots,N\}$ are continuous in $\lambda$.

\section{Proof of Lemma \ref{lem:sublipschitz}}\label{app:2}
We begin with a technical lemma. Let $Q_c: = [Q_1,Q_2,\dots,Q_N]^{\top}, R_c: = [R_1,R_2,\dots,R_N]^{\top}$, and $G_{\lambda} = BR_{\lambda}^{-1}B^{\top}$.
\begin{lemma}\label{lem:lip}
	The following results hold.
	
	(a) For $\lambda \geq 0$, it holds that $\|A-BK_{\lambda}^*\| \leq \|A\|$ and $\|K_{\lambda}^*\| \leq  D(\lambda)\|A\|\|B\|/{\underline{\sigma}(R_0)} $.
	
	(b) For any $\lambda, \lambda'\in \Omega$, if $\|\lambda' - \lambda\| \leq \underline{\sigma}(R_0)/(2\|R_c\|)$, it holds that
	$\|G_{\lambda'}- G_{\lambda} \| \leq ({2\|R_c\|\|B\|^2}/\underline{\sigma}(R_0))\|\lambda' - \lambda \|.$
	
	(c) For any $\lambda, \lambda' \geq 0$, it holds  $\|Q_{\lambda'} -Q_{\lambda }\| \leq \|Q_c\|\|\lambda' - \lambda\|$.
\end{lemma}
\begin{proof}
	(a) By the definition of $K_{\lambda}^*$ in \eqref{def:K} and matrix inverse lemma, it holds that
	\begin{align*}
	\|A-BK_{\lambda}^*\|
	&= \|A - B(R_0+\lambda R_c + B^{\top}P_{\lambda}^*B)^{-1}B^{\top}P_{\lambda}^*A \| \\
	& = \|(I + B(R_0 + \lambda R_c)^{-1}B^{\top}P_{\lambda}^*)^{-1}A \| 
	\leq \|A\|, \\
	\| K_{\lambda}^* \| &= \|(R_0+\lambda R_c + B^{\top}P_{\lambda}^*B)^{-1}B^{\top}P_{\lambda}^*A \|\\
	& \leq \frac{1}{\underline{\sigma}(R_0)}\|A\|\|B\|\|P_{\lambda}^*\| \leq \frac{\mathcal{L}_{\lambda}^*}{\underline{\sigma}(R_0)}\|A\|\|B\|.
	\end{align*}
	
	(b) By the perturbation theorem of matrix inverse~\cite{stewart1990matrix}, if $\|\lambda'  - \lambda\| \leq \underline{\sigma}(R_0)/(2\|R_c\|)$, then it holds that
	\begin{align*}
	\|G_{\lambda'}- G_{\lambda} &\| \leq \|B\|^2 \|(R_0 + \lambda'  R_c)^{-1} - (R_0+\lambda R_c)^{-1}  \| \\
	&\leq \frac{2\|R_c\|\|B\|^2}{\underline{\sigma}(R_0)} \|\lambda'  - \lambda \|.
	\end{align*}
	
	(c)  
	$
	\|Q_{\lambda'} -Q_{\lambda }\| = \| \lambda'  Q_c -  \lambda Q_c\| \leq \|Q_c\|\|\lambda'  - \lambda\|.
	$
\end{proof}

\begin{lemma}\label{lem:P}
	For any $\lambda, \lambda'\in \Omega$, there exists some constant $\mathcal{B}_{\Omega}^{\lambda}>0$ such that if $\|\lambda'-\lambda\| \leq \mathcal{B}_{\Omega}^{\lambda}$, it holds that 
	\begin{equation}
	\|P_{\lambda'}^*-P_{\lambda}^*\| \leq \mathcal{B}_{\Omega}^{P} \|\lambda' - \lambda\|,
	\end{equation}
	for some constant $ \mathcal{B}_{\Omega}^{P}>0$.
\end{lemma}
\begin{proof}
	Since the ARE (\ref{equ:dare}) has a unique positive definite solution for any $\lambda \geq 0$, we use perturbation theory for Riccati equations \cite[Theorem 4.1]{sun1998perturbation} for the proof. For brevity, define $\Delta G = G_{\lambda'}-G_{\lambda}, \Delta Q= Q_{\lambda'}-Q_{\lambda}$ and the following quantities
	\begin{equation}\label{def:quant}
	\begin{aligned}
	&g=\left\|G_{\lambda}\right\|, ~H=P_{\lambda}^*\left({I}+G_{\lambda} P_{\lambda}^*\right)^{-1} A ,~ \delta=\frac{H\|\Delta G\|}{1-\psi\|\Delta G\|}, \\ &f=\left\|\left({I}+G_{\lambda} P_{\lambda}^*\right)^{-1}\right\|,   ~
	\phi=\left\|\left({I}+G_{\lambda} P_{\lambda}^*\right)^{-1} A\right\|, \\
	& \beta=f \delta(2 \phi+f \delta),~\psi=\left\|P_{\lambda}^* \cdot\left({I}+G_{\lambda} P_{\lambda}^*\right)^{-1}\right\| \\ 
	 &
	T_{\lambda}={I}-\left(A-B K_{\lambda}^*\right)^{\top} \otimes\left(A-B K_{\lambda}^*\right)^{\top}, \\
	&l=\left\|T_{\lambda}^{-1}\right\|^{-1},  
	q= \left\|T_{\lambda}^{-1}\left(H^{\top} \otimes H\right)\right\|,\theta=\frac{l}{\phi+\sqrt{\phi^2+l}},\\
	& \kappa=\frac{1}{l}\left\|\Delta Q\right\|+\left(q+\frac{H \psi \delta}{l}\right)\|\Delta G\|, ~
	\alpha=\frac{f\left\|A\right\|}{1-\psi\|\Delta G\|},\\
	& p = \|Q_c\|{\underline{\sigma}(R_0)}+2\left(ql+H \psi \delta\right){\|R_c\|\|B\|^2}.
	\end{aligned}
	\end{equation}
	
	It follows from the definition that $l$ is finite and 
	$
	l = 1/\|T^{-1}_{\lambda}\| = \underline{\sigma}(T_{\lambda}) > 0.
	$
	Moreover, since $\underline{\sigma}(T_{\lambda})$ is continuous in $\lambda$, $l$ is uniformly bounded above zero over the compact set $\Omega$. 	If we require the following conditions 
	\begin{equation}\label{equ:b19}
	\|\Delta G\| \leq \frac{1}{\psi}, \quad 1-f g \xi_* \geq 0, \quad \frac{f \delta+\phi f g \xi_*}{1-f g \xi_*} \leq \theta, 
	\end{equation}
	where $\xi_*=(2 l \kappa) \cdot(l / 2+l f g \kappa)^{-1}$, then the condition (4.40) of \cite{sun1998perturbation} is satisfied. Additionally, if we require that 
	\begin{equation}\label{equ:b20}
	\beta  \leq \frac{fH\|\Delta G\|_2}{1-\psi\|\Delta G\|_2}(2\phi + 2l + \frac{fH\|\Delta G\|_2}{1-\psi\|\Delta G\|_2})  \leq  l/2,
	\end{equation}
	then the definition of $\xi_*$ here is strictly larger than that in \cite{sun1998perturbation}. If we also let the following holds
	\begin{equation}\label{equ:b21}
	\kappa= \frac{1}{l}\left\|\Delta Q\right\|+\left(q+\frac{H \psi \delta}{l}\right)\|\Delta G\| \leq \frac{(l / 2)^2}{2 l f g(l+2 \alpha)},
	\end{equation}
	then the condition in (4.41) of \cite{sun1998perturbation} holds, since \eqref{equ:b20}-\eqref{equ:b21} implies 
	$$
	\kappa  \leq 
	\frac{(l-\beta)^2}{l f g\left(l-\beta+2 \alpha+\sqrt{(l-\beta+2 \alpha)^2-(l-\beta)^2}\right)}.
	$$
	Under \eqref{equ:b19}-\eqref{equ:b21},  \cite[Theorem 4.1]{sun1998perturbation} implies that for $\|\lambda'-\lambda\| \leq \underline{\sigma}(R_0)/(2\|R_c\|)$, we have
	\begin{equation}\label{equ:p}
	\begin{aligned}
	&\|P_{\lambda'}^* - P_{\lambda}^*\| \\
	&\leq \frac{2 l \kappa}{l / 2+l f g \kappa} \leq 4\kappa = \frac{4}{l}\left\|\Delta Q\right\|+\left(4q+\frac{4H \psi \delta}{l}\right)\|\Delta G\|\\ 
	& \leq \left(\frac{4\|Q_c\|}{l} +  4q+\frac{8H \psi \delta {\|R_c\|\|B\|^2}}{l\underline{\sigma}(R_0)} \right)\|\lambda' - \lambda \|.
	\end{aligned}
	\end{equation}
	
	Now, we provide sufficient conditions of \eqref{equ:b19}-\eqref{equ:b21} such that the bound in (\ref{equ:p}) holds. Note that (\ref{equ:b19}) is equivalent to 
	\begin{equation}\label{equ:23}
	f g \kappa \leq 1 / 2, ~\text{and} ~~f \delta+(\phi+\theta) f g \xi_* \leq \theta.
	\end{equation}
	Since $\xi_* \leq 4 \kappa$, $f\delta\geq 0$, and $fg\kappa \geq 0$, if we let 
	\begin{equation}\label{equ:24}
	f \delta+4(\phi+\theta) f g \kappa \leq \theta,
	\end{equation}
	then (\ref{equ:23}) holds. Thus, we only need the following sufficient condition for (\ref{equ:24}) 
	$$
	\frac{fH\|\Delta G\|}{1-\psi\|\Delta G\|} + \frac{4(\phi + \theta)fg}{l}\left(  \|\Delta Q \|+ (ql+ {H \psi \delta}  )\|\Delta G\|\right) \leq \theta.
	$$
	
	By the second statement of Lemma \ref{lem:lip}, it suffices to let 
	\begin{align}
	 &({2\|R_c\|\|B\|^2}/{\underline{\sigma}(R_0)})\cdot \|\lambda'-\lambda\| \leq  1/(2\psi) \label{equ:21}\\ 
	 &2fH\|\Delta G\| +  {4(\phi + \theta)fg}\left(  \|\Delta Q \|+ (ql+ {H \psi \delta}  )\|\Delta G\|\right)/l \leq \theta. \notag
	\end{align}

	Using the bounds in Lemma \ref{lem:lip} for $\|\Delta G\|$ and $\|\Delta Q \|$, if
	\begin{equation}\label{equ:16}
	\begin{aligned}
	&\|\lambda'-\lambda\| \leq 
	\min\bigg\{ \frac{\theta l{\underline{\sigma}(R_0)}}{4fH{\|R_c\|\|B\|^2}l + 4(\phi + \theta)fgp},\\
	&\frac{\underline{\sigma}(R_0)}{4\psi {\|R_c\|\|B\|^2}}, \frac{\underline{\sigma}(R_0)}{2\|R_c\|}\bigg\},
	\end{aligned}
	\end{equation}
	then (\ref{equ:b19}) holds.
	
	To ensure (\ref{equ:b20}), it suffices to let \eqref{equ:21} hold and
	\begin{equation}
	{2fH\|\Delta G\|_2}(2\phi + 2l + {2fH\|\Delta G\|_2})  \leq  l/2,
	\end{equation}
	which holds for
	\begin{equation}\label{equ:17}
	\|\lambda'-\lambda\| \leq \min \left\{\frac{\underline{\sigma}(R_0)}{4\psi {\|R_c\|\|B\|^2}}, \frac{\sqrt{l+2(\phi+l)^2}}{\sqrt{8}fH} - \frac{\phi + l}{2fH}\right\} .
	\end{equation}
	
	Since $\alpha \geq f\|A\|$, the condition (\ref{equ:b21}) holds by letting
	\begin{equation}\label{equ:27}
	\frac{1}{l}\left\|\Delta Q\right\|+\left(q+\frac{H \psi \delta}{l}\right)\|\Delta G\| \leq \frac{l}{8fg(l+2f\|A\|)}.
	\end{equation}
	By Lemma \ref{lem:lip}, (\ref{equ:27}) holds if
	\begin{equation}\label{equ:18}
	\|\lambda'-\lambda\| \leq \frac{l^2\underline{\sigma}(R_0)}{8fgp(l+2f\|A\|) }.
	\end{equation} 
	
	Thus, under (\ref{equ:16}), (\ref{equ:17}) and (\ref{equ:18}), the bound (\ref{equ:p}) holds. Finally, we note that upper bounds for $\|\lambda'-\lambda\|$ in (\ref{equ:16}), (\ref{equ:17}) and (\ref{equ:18}) are all lower bounded above zero, since (a) all the quantities in (\ref{def:quant}) are positive, and (b) they are continuous in $\lambda$ and hence have uniform bounds over the compact set $\Omega$. Thus, it suffices to let $\mathcal{B}_{\Omega}^{\lambda}$ be the lower bounds of (\ref{equ:16}), (\ref{equ:17}) and (\ref{equ:18}) over $\Omega$. By (\ref{equ:p}), we let
	$$
	\mathcal{B}_{\Omega}^P := \frac{4\|Q_c\|}{l} +  4q+\frac{8H \psi \delta {\|R_c\|\|B\|^2}}{l\underline{\sigma}(R_0)}
	$$
	such that $\|P_{\lambda'}^* - P_{\lambda}^*\| \leq \mathcal{B}_{\Omega}^P \|\lambda' -\lambda\|$.
\end{proof}

\begin{lemma}\label{lem:lip_K}
	For any $\lambda, \lambda'\in \Omega$, if $\|\lambda'-\lambda\| \leq \mathcal{B}_{\Omega}^{\lambda}$, then we have  
	$
	\|K_{\lambda'}^* - K_{\lambda}^*\| \leq  \mathcal{B}_{\Omega}^{K}  \|\lambda' -\lambda \|
	$
	for some constant $\mathcal{B}_{\Omega}^{K}$.
\end{lemma}
\begin{proof}
	By the definition of $K_{\lambda}^*$, it holds that
	\begin{align*}
	&(R_0 + \lambda R_c+B^{\top}P_{\lambda}^*B)K_{\lambda}^*= B^{\top}P_{\lambda}^*A\\
	&(R_0 + \lambda'  R_c+B^{\top}P_{\lambda'}^*B)K_{\lambda'}^*= B^{\top}P_{\lambda'}^*A.
	\end{align*}
	
	Subtracting the first equation from the second yields that
	\begin{align*}
	&\left(R_0 + \lambda' R_c+B^{\top}P_{\lambda'}^*B\right)\left(K_{\lambda'}^*-K_{\lambda}^*\right)= \\ & B^{\top}(P_{\lambda'}^*-P_{\lambda}^*)A-\left(R_c(\lambda' -\lambda) + B^{\top}(P_{\lambda'}^*-P_{\lambda}^*)B\right)K_{\lambda}^*.
	\end{align*}
	
	By Lemmas \ref{lem:lip} and \ref{lem:P},  for any $\lambda,\lambda' \in \Omega$ and $\|\lambda'-\lambda\| \leq \mathcal{B}_{\Omega}^{\lambda}$, it holds that
	\begin{align*}
	&\underline{\sigma}(R_0)\|K_{\lambda'}^*-K_{\lambda}^*\| \\
	&\leq  \left\|  B^{\top}(P_{\lambda'}^*-P_{\lambda}^*)(A-BK_{\lambda}^*)- (\lambda' -\lambda)R_cK_{\lambda}^* \right\| \\
	&\leq  \|B\|\|A-BK_{\lambda}^*\|\|P_{\lambda'}^*-P_{\lambda}^*\| + \|R_c\|\|K_{\lambda}^*\| \|\lambda' -\lambda\|\\
	&\leq   \|A\|\|B\|\|P_{\lambda'}^*-P_{\lambda}^*\| + \frac{D(\lambda)}{\underline{\sigma}(R_0)}\|A\|\|B\|\|R_c\| \|\lambda' -\lambda\|   \\
	&\leq  \left(\|A\|\|B\|\mathcal{B}_{\Omega}^{P} + \frac{D^*}{\underline{\sigma}(R_0)}\|A\|\|B\|\|R_c\| \right) \|\lambda' -\lambda\|.
	\end{align*}

	Then, letting $$\mathcal{B}_{\Omega}^{K} = \frac{1}{\underline{\sigma}(R_0)}\left(\|A\|\|B\|\mathcal{B}_{\Omega}^{P} + \frac{D^*}{\underline{\sigma}(R_0)}\|A\|\|B\|\|R_c\| \right)$$
	completes the proof.
\end{proof}

%
%\begin{lemma} 
%	For any $\lambda, \lambda'\in \Omega$, there exist constants $\gamma>0$ and $\mu >0$ such that if $\|\lambda'-\lambda\| \leq \gamma$, 
%	then it holds that 
%	$
%	\|\nabla D(\lambda') - \nabla D(\lambda)\| \leq \mu \|\lambda' -{\lambda}\|.
%	$
%\end{lemma}  
	Let $P_i(K)$ be the solution to the Lyapunov equation
	$
	P_i(K) = Q_i + K^{\top}R_iK + (A-BK)^{\top}P_i(K)(A-BK).
	$
	By the perturbation theory of Lyapunov equations \cite[Lemma 27]{fazel2018global}, if
	\begin{equation}\label{equ:KK}
	\|K_{\lambda'}^* - K_{\lambda}^*\| \leq \frac{\underline{\sigma}(Q_0)}{4J_i(K_{\lambda}^*)\|B\|(\|A\|+1)},
	\end{equation}
	then it holds that
	\begin{equation}\label{equ:31}
	\begin{aligned}
	&|J_i(K_{\lambda'}^*) - J_i(K_{\lambda}^*)| \leq\|P_i(K_{\lambda'}^*) - P_i(K_{\lambda}^*)\| \\
	&\leq 6\|K_{\lambda}^*\|\|R_i\|\left(\frac{J_i(K_{\lambda}^*)}{\underline{\sigma}(Q_0)} \right)^2\|K_{\lambda'}^*-K_{\lambda}^*\| \\
	&\times(\|K_{\lambda}^*\|\|B\|\|A-BK_{\lambda}^*\| + \|K_{\lambda}^*\|\|B\|+1)\\
	&\leq 6\|R_i\|\left(\frac{D(\lambda)\|A\|\|B\|}{\underline{\sigma}(R_0)}\right)\left(\frac{J_i(K_{\lambda}^*)}{\underline{\sigma}(Q_0)} \right)^2\|K_{\lambda'}^*-K_{\lambda}^*\| \\
	&\times \left(\left(\frac{D(\lambda)\|A\|\|B\|^2}{\underline{\sigma}(R_0)}\right) (\|A\|+1)+1 \right).
	\end{aligned}
	\end{equation} 
	
	Since $J_i(K_{\lambda}^*)$ is continuous in $\lambda$ over the compact set $\Omega$, it is uniformly upper bounded by some constant $\mathcal{B}_{\Omega}^{J_i}$, i.e., $J_i(K_{\lambda}^*) \leq \mathcal{B}_{\Omega}^{J_i}$. Noting that $D(\lambda)\leq D^*$, \eqref{equ:31} can be further bounded as
	$
	|J_i(K_{\lambda'}^*) - J_i(K_{\lambda}^*)| \leq \mathcal{B}_{\Omega}^{i}\|K_{\lambda'}^*-K_{\lambda}^*\| 
	$,
	where 
	\begin{equation}\label{equ:def_Bi}
	\begin{aligned}
	\mathcal{B}_{\Omega}^{i} &= 6\|R_i\|\left(\frac{D^*\|A\|\|B\|}{\underline{\sigma}(R_0)}\right)\left(\frac{\mathcal{B}_{\Omega}^{J_i}}{\underline{\sigma}(Q_0)} \right)^2 \\
	& \times \left(\left(\frac{D^*\|A\|\|B\|^2}{\underline{\sigma}(R_0)}\right) (\|A\|+1)+1 \right).
	\end{aligned}
	\end{equation}

	By Lemma \ref{lem:lip_K}, a sufficient condition for \eqref{equ:KK} is 
	$$
	\|\lambda'-\lambda\| \leq \min \left\{ \mathcal{B}_{\Omega}^{\lambda},   \frac{\underline{\sigma}(Q_0)}{4\mathcal{B}_{\Omega}^{K}   \max\{\mathcal{B}_{\Omega}^{J_i}\} \|B\|(\|A\|+1)} \right\}.
	$$

	Then, it follows from the definition of the subgradient $\nabla D(\lambda)$ and Lemma \ref{lem:lip_K} that
	$
	\|\nabla D(\lambda') - \nabla D(\lambda)\| \leq \sum_{i=1}^{N}|J_i(K_{\lambda'}^*) - J_i(K_{\lambda}^*)|  
	\leq  \sum_{i=1}^{N}\mathcal{B}_{\Omega}^{i}\|K_{\lambda'}^*-K_{\lambda}^*\| \leq \sum_{i=1}^{N}\mathcal{B}_{\Omega}^{i}\mathcal{B}_{\Omega}^{K}  \|\lambda' -\lambda \|
	$.  
	
	\section{Proof in Section \ref{subsec:conver}} \label{app:3}
	\subsection{Proof of Lemma \ref{lem:sub_error}}
		The proof follows from the derivation in \eqref{equ:KK}-\eqref{equ:31}. Let 
		$$
		\gamma_2 = \frac{\underline{\sigma}(Q_0)}{4 \max_i\{\mathcal{B}_{\Omega}^{J_i}\}\|B\|(\|A\|+1)}.
		$$
		If  $\epsilon \leq \gamma_2$, then $\|d^k - \nabla D(\lambda^k)\| \leq\sum_{i=1}^{N}\mathcal{B}_{\Omega}^{i} \|K^k-K^*_{\lambda^k} \|\leq \sum_{i=1}^{N}\mathcal{B}_{\Omega}^{i}  \epsilon:= c\epsilon$, where $\mathcal{B}_{\Omega}^{i}$ is defined in \eqref{equ:def_Bi}.

		Since $\|\nabla D(\lambda)\|$ is continuous in $\lambda$ and $\Omega$ is a compact set, it follows that $\|\nabla D(\lambda^k)\|$ and $\|d^k\|$ are uniformly upper bounded.

\subsection{Proof of Theorem \ref{thm:finald}}

\begin{lemma}\label{lem:boudn}
Let $\lambda,y \in \Omega$, $\lambda^+ = \Pi_{\Omega}(\lambda + \eta d_{\lambda})$,  $\|d_{\lambda}-\nabla D(\lambda)\| \leq c\epsilon$, $d_{\lambda} \leq \bar{d}$, $g_{\lambda}= (\lambda^+ - \lambda)/\eta$. Then, if $\eta \leq \gamma_1/\bar{d}$, we have
$$
D(\lambda^+) - D(y) \geq g_{\lambda}^{\top}(\lambda-y) - c\epsilon\omega + (\eta - \frac{\mu\eta^2}{2})\|g_{\lambda}\|^2.
$$
\end{lemma}
\begin{proof}
	By the non-expansiveness of the projection, we have
	$
	(\lambda^+ - y)^{\top}(\lambda^+-\lambda -\eta d_{\lambda})<0,
	$
	which is equivalent to
	$
	d_{\lambda}^{\top}(\lambda^+-y)\geq g_{\lambda}^{\top}(\lambda^+-y).
	$
	Then, it follows that
	\begin{align*}
	&D(\lambda^+) - D(y) = D(\lambda^+) - D(\lambda) + D(\lambda) - D(y) \\
	& \geq  \nabla D(\lambda)^{\top}(\lambda^+-\lambda) - \frac{\mu}{2}\|\lambda^+-\lambda\|^2 + \nabla D(\lambda)^{\top} (\lambda - y) \\
	&=  \nabla D(\lambda)^{\top}(\lambda^+ -y)  - \frac{\mu\eta^2}{2}\|g_{\lambda}\|^2 \\
	& =
	(d_{\lambda}+\nabla D(\lambda)-d_{\lambda})^{\top}(\lambda^+-y) - \frac{\mu\eta^2}{2}\|g_{\lambda}\|^2 \\
	&\geq g_{\lambda}^{\top}(\lambda^+-y) - \|\nabla D(\lambda)-d_{\lambda}\|\|\lambda^+-y\|  - \frac{\mu\eta^2}{2}\|g_{\lambda}\|^2 \\
	&\geq  g_{\lambda}^{\top}(\lambda-y) - c\epsilon\omega + (\eta - \frac{\mu\eta^2}{2})\|g_{\lambda}\|^2,
	\end{align*}
	where the first inequality follows from Lemma \ref{lem:sublipschitz} and the concavity of $D(\lambda)$.
\end{proof}

By Lemma \ref{lem:boudn}, it holds that
\begin{align}
&D(\lambda^{k+1}) - D(\lambda^k) \geq  - c\epsilon\omega + (1 - \frac{\mu\eta}{2})\|\lambda^{k+1}-\lambda^k\|^2 \label{subequ:1} \\
&D(\lambda^{k+1}) - D(\lambda^*)  \geq -\frac{1}{\eta}\|\lambda^{k+1}-\lambda^k\|\omega -c\epsilon \omega.\label{subequ:2}
\end{align}

Rearranging \eqref{subequ:1} yields
\begin{align*}
&D(\lambda^{k+1}) - D(\lambda^k) \geq (1-\frac{\mu \eta}{2})( \|\lambda^{k+1}-\lambda^k\|^2 + \eta c\epsilon )^2 \\
& -(2-\mu \eta)\eta c\epsilon \|\lambda^{k+1}-\lambda^k\| - (1-\frac{\mu \eta}{2})(\eta c \epsilon)^2 - c\epsilon \omega\\
&\geq (1-\frac{\mu \eta}{2})( \|\lambda^{k+1}-\lambda^k\|^2 + \eta c\epsilon )^2  - p_1.
\end{align*}

Combining \eqref{subequ:2}, it leads to that
$$
D(\lambda^{k+1}) - D(\lambda^k) \geq \frac{1}{p_2}(D^*-D(\lambda^{k+1}))^2 - p_1.
$$
 Rearranging it yields
$$
(D^*-D(\lambda^{k+1}))^2 \leq p_2 (D(\lambda^{k+1}) - D(\lambda^k)) + p_1p_2.
$$
Summing up both sides from $0$ to $k$ yields
$$
\sum_{i=0}^{k-1}(D^*-D(\lambda^{i+1}))^2 \leq p_2 D(\lambda^{k}) + kp_1p_2.
$$
Dividing both sides by $k$ and using the AM-GM inequality yield that
$$
\frac{1}{k}\sum_{i=0}^{k-1}(D^*-D(\lambda^{i+1})) \leq \frac{\sqrt{p_2 D^*}}{\sqrt{k}} + \sqrt{p_1p_2}.
$$

\bibliographystyle{IEEEtran}
\bibliography{mybibfile}
\end{document}